\documentclass[12pt,a4paper,reqno]{amsart}
\usepackage{a4wide,amsmath,amscd,amssymb,latexsym}

\def\q{\hskip0.17cm}
\def\,{\hskip0.12cm}

\newtheorem{Thm}{Theorem}[section]

\newtheorem{Lem}[Thm]{Lemma}

\theoremstyle{definition}

\newtheorem{Exm}[Thm]{Example}
\theoremstyle{remark}
\newtheorem{Rmk}[Thm]{Remark}

\newcommand{\sm}{\left(\smallmatrix}
\newcommand{\esm}{\endsmallmatrix\right)}

\begin{document}
 \title[Borcherds products associated with certain Thompson series]{Borcherds
 products associated with certain Thompson series}
 \author{Chang Heon Kim}
 \address{
 Max-Planck-Institut f\"ur Mathematik, Vivatsgasse 7, 53111 Bonn,
 Germany}
 \email{chkim@mpim-bonn.mpg.de}
 \maketitle
 \renewcommand{\thefootnote}%
             {}

 \begin{abstract}
 We apply Zagier's result for the traces of singular moduli to
 construct Borcherds products in higher level cases.
 \end{abstract}
 \section{Introduction}
 \par
 Let $M_{1/2}^!$ be the additive group consisting of nearly
 holomorphic modular forms of weight $1/2$ for $\Gamma_0(4)$ whose
 Fourier coefficients are integers and satisfy
 the Kohnen's ``plus space" condition
 (i.e. $n$-th coefficients vanish
 unless $n\equiv 0$ or $1$
 modulo 4).
  We also let $\mathcal B$ be the multiplicative group
 consisting of meromorphic modular forms for some characters of
 $SL_2(\mathbb Z)$ of integral weight with leading coefficient 1
 whose coefficients are integers and all of whose zeros and poles
 are either cusps or imaginary quadratic irrationals. Borcherds
 \cite{Borcherds} gave an isomorphism between $M_{1/2}^!$ and $\mathcal B$
 \, by means of infinite products
 which we call modular products or Borcherds products.
 \par
 Let $d$ denote a positive integer congruent to 0 or 3 modulo 4.
 We denote by ${\mathcal Q}_d$ the set of positive definite binary quadratic
 forms $Q=[a,b,c]=aX^2+bXY+cY^2 \, (a,b,c\in \mathbb Z)$ of
 discriminant $-d$, with usual action of the modular group
 $\Gamma=PSL_2(\mathbb Z)$. To each $Q\in {\mathcal Q}_d$, we associate its
 unique root $\alpha_Q\in {\mathfrak H}$ (=upper half plane). We
 define the {\it Hurwitz-Kronecker class number} $H(d)$ by
 $H(d)=\underset{Q\in {\mathcal Q}_d/\Gamma}{\sum} \frac{1}{w_Q}$ where
 $w_Q=|\Gamma_Q|$.
 For instance, we have $H(3)=1/3$, $H(4)=1/2$,
 $H(7)=H(8)=H(11)=1$, $H(12)=4/3$, $H(15)=2$, etc.
 For the modular invariant $j(\tau)$,
 we define a function $\mathcal
 H_d(j(\tau))\in \mathcal B$ \,
 by $\underset{Q\in {\mathcal Q}_d/\Gamma}{\prod}
 (j(\tau)-j(\alpha_Q))^{1/w_Q}$.
 On the other hand, for each $d$ there is a unique modular form
 $f_{d,1}\in M_{1/2}^!$ having a Fourier development of the form
 $ f_{d,1}=q^{-d}+\sum_{D>0} A(D,d) q^D, \, q=e^{2\pi i\tau}
  (\tau\in \mathfrak H). $ Then Borcherds'
 theorem says that
 \begin{equation}
 \mathcal
 H_d(j(\tau))=q^{-H(d)} \prod_{u=1}^\infty (1-q^u)^{A(u^2,d)}
 \tag{$*$}
 \end{equation}
 Zagier \cite{Zagier} described the trace of a singular modulus of
 discriminant $-d$ ($=\underset{Q\in {\mathcal Q}_d/\Gamma}{\sum} \frac{1}{w_Q}
  (j(\alpha_Q)-744)$)
  as the coefficient of $q^d$ in a fixed
 modular form $-g_{1,1}(\tau)$ of weight $3/2$. By making use of
 this formula and considering Hecke operators in integral and
 half-integral weight Zagier reproved $(*)$ (see \cite{Zagier} \S 6).
 Moreover he generalized the trace formula to the group
 $\Gamma_0(N)^*$ (=the group generated by $\Gamma_0(N)$ and all
 Atkin-Lehner involutions $W_e$ for $e||N$) for $2\le N \le 6$
 (see \cite{Zagier} \S 8).
 \par
 In this article we find an analogue of $(*)$
 in higher level cases $N=2,3,5,6$
 by applying
 Zagier's Theorem 8 in \cite{Zagier}.
 Let $M_{k-1/2}^{+\cdots+}(N)^!$ be the vector space consisting
 of nearly holomorphic modular forms of half-integral weight
 $k-1/2$ on $\Gamma_0(4N)$ whose $n$-th Fourier coefficient
 vanishes unless $(-1)^{k-1}n$ is a square modulo $4N$. There is
 a unique modular form $f_{d,N}\in M_{1/2}^{+\cdots+}(N)^!$
 having a Fourier expansion of the form
 $$ f_{d,N}=q^{-d}+\sum_{D>0} A(D,d) q^D. $$
 An explicit construction of $f_{d,N}$ is given in the appendix
 and the uniquess of $f_{d,N}$ is shown in the end of \S 2.
 Let ${\mathcal Q}_{d,N}$ be the set of
 forms $Q=[a,b,c]\in {\mathcal Q}_d$ satisfying $N|a$. Then $\Gamma_0(N)^*$
 naturally acts on ${\mathcal Q}_{d,N}$ and the
 quotient ${\mathcal Q}_{d,N}/\Gamma_0(N)^*$
 has a bijection with ${\mathcal Q}_d/\Gamma$ (see \cite{Zagier} \S 8).
 We can therefore define, for the Hauptmodul $t(\tau)$ for
 $\Gamma_0(N)^*$, a modular function
 $\mathcal
 H_d(t(\tau))$ by $ \underset{Q\in {\mathcal Q}_{d,N}/\Gamma_0(N)^*}{\prod}
 (t(\tau)-t(\alpha_Q))^{1/w_Q}$.
 In \S 3 we will prove the following theorem.

 \begin{Thm} Let $1\le N \le 6$ other than 4 and $t$ be the
 Hauptmodul for $\Gamma_0(N)^*$. Let $-d$ be the discriminant
 corresponding to a Heegner point (i.e. the discriminant of $Q\in
 {\mathcal Q}_{d,N}$ with the condition that if $f^2$ divides $d$, then
 $(f,N)=1$).
  Define $A^*(u^2,d)=2^{s(u,N)} A(u^2,d)$ where $s(u,N)$ is
 the number of distinct prime factors dividing $(u,N)$. Then
 $$ {\mathcal H}_d(t(\tau))=q^{-H(d)}\prod_{u=1}^\infty
 (1-q^u)^{A^*(u^2,d)}.$$
 \label{Borcherds}
 \end{Thm}
 We remark that this theorem is related to the problem of
 generalizing Borcherds'
 theorem (\cite{Borcherds} Theorem 14.1) to higher levels
 (\cite{Borcherds}
 problem 10 in \S 17).
 In some sense Borcherds proved it himself in \cite{Borcherds2}
 Theorem 13.3. The vector valued modular forms he uses include the
 higher level case because a higher level form can be induced upto
 a vector valued form of level $1$. An explicit infinite product
 is given in part 5 of Theorem 13.3 of \cite{Borcherds2}. But as
 he pointed out, it seems to take a bit of effort to unravel it to
 see what it says in the case of modular forms.
 Also Bruinier \cite{Bruinier} proved that every automorphic forms
 with zeros on Heegner divisors can be written as modular
 products in the case that the lattice considered splits two hyperbolic
 planes over $\mathbb Z$.
 \par
 Finally in \S 4 by using the idea given in \cite{KKKO}
 we derive a recursion formula which enables us to
 estimate all $A^*(u^2,d)$ for $u\ge 1$ from the Fourier
 coefficients of ${\mathcal H}_d(t(\tau))$.

\section{Preliminaries}
 \subsection{Generalized Hecke operator}
 \par
 Let $N$ be a positive integer and $e$ be any Hall divisor of
$N$ (written $e||N$), that is, a positive divisor of $N$ for which
$(e, N/e)=1$. We denote by $W_{e,N}$ a matrix $\sm ae & b \\ cN &
de \esm$ with $\det W_{e,N}=e$ and $a,b,c,d\in \mathbb Z$, and
call it an {\it Atkin-Lehner involution}. Let $S$ be a subset of
Hall divisors of $N$ and let $N+S$ be the subgroup of
$PSL_2(\mathbb R)$ generated by $\Gamma_0(N)$ and all Atkin-Lehner
involutions $W_{e,N}$ for $e\in S$
 (we may choose $S$ so that $1\not\in S$ and if
 $e_1, e_2 \in S$, then $e_1 e_2/(e_1, e_2)^2\in S$
 unless $e_1=e_2$). We assume that the genus of
the group $N+S$ is zero. Then there exists a unique modular
function $t$ with respect to $N+S$ satisfying
\par\noindent (i) $t$ is holomorphic
 on the complex upper half plane $\mathfrak H$,
\par\noindent (ii) $t$ has the Fourier expansion at $\infty$ of
the form $$ t= q^{-1}+\sum_{k\ge 1} H_k q^k, \, q=e^{2\pi i \tau}
\, (\tau\in \mathfrak H), $$
\par\noindent (iii) $t$ is holomorphic at all cusps which are not
equivalent to $\infty$ under $N+S$. Such a function $t$ is called
the {\it Hauptmodul} \, for $N+S$.
 By the result of Borcherds
 \cite{Borcherds92} $t$ becomes a monstrous function whose
 Fourier coefficients are related to representations of the monster group
$\mathbb M$ except for the three cases ($25-$, $49+49$ and
$50+50$). More precisely the $q$-series of $t$ coincides with a
Thompson series $T_g(\tau)=\sum_{n\in \mathbb Z} \text{Tr}(g|V_n)
q^n$ for some element $g$ of $\mathbb M$. Here $V=\bigoplus_{n\in
\mathbb Z} V_n$ is the infinite dimensional graded representation
of $\mathbb M$ constructed by Frenkel {\it et al.}
\cite{Frenkel1, Frenkel2}.
 For a prime number $p$, let $t^{(p)}$ be the Hauptmodul for
 $N^{(p)}+S^{(p)}$ where $N^{(p)}=N/(p,N)$ and $S^{(p)}$ is the
 set of all $e$ in $S$ which divide $N^{(p)}$. Generally if $m=p_1
 p_2\cdots p_r$ is a product of primes $p_i$, then we define the
 {\it $m$-th replicate} $t^{(m)}$ of $t$ by
 $$t^{(m)}=(\cdots ((t^{(p_1)})^{(p_2)}\cdots)^{(p_r)}.$$
 For every positive integer $n$, let $t_n$ be a unique polynomial
 of $t$ satisfying $t_n\equiv q^{-n} \mod q \mathbb C[[q]]$.
 Define
 the {\it $m$-th generalized Hecke operator} $T(m)$
 \cite{ACMS, Ferenbaugh, Fer2, Koike} by
 $$ t_n|_{T(m)}=\sum_{ad=m\atop 0\le b < d} t_n^{(a)}\left(
 \frac{a\tau +b}{d}\right).$$
 The $m$-th replication formula \cite{Fer2, Koike} says that $t_m=t|_{T(m)}$.
 \subsection{Jacobi forms}
 \par
 A {\it (holomorphic) Jacobi form on $SL_2(\mathbb Z)$} is defined to be a
 holomorphic function $\phi : \mathfrak H \times \mathbb C \to
 \mathbb C$ \, satisfying the two transformation equations
 \begin{align*}
 \phi \left(\frac{a\tau+b}{c\tau+d}, \frac{z}{c\tau+d}\right)
 & = (c\tau+d)^k e^{\frac{2\pi i Ncz^2}{c\tau+d}} \phi (\tau,z)
     \q \q (\sm a&b \\c&d \esm \in SL_2(\mathbb Z)), \\
 \phi (\tau, z+\lambda\tau+\mu)
 & = e^{-2\pi iN(\lambda^2 \tau+2\lambda z)} \phi (\tau,z)
     \q \q (\sm \lambda & \mu \esm \in \mathbb Z^2)
 \end{align*}
 and having a Fourier expansion of the form
 \begin{equation}
 \phi (\tau,z)=\sum_{n,r\in \mathbb Z \atop 4Nn-r^2\ge 0} c(n,r) q^n
 \zeta^r \q \q \q (q=e^{2\pi i\tau}, \, \zeta=e^{2\pi iz}).
 \label{Fourier}
 \end{equation}
 Here $k$ and $N$ are positive integers, called the {\it weight}
 and {\it index} of $\phi$, respectively.
 The coefficient $c(n,r)$ depends only on
 $4Nn-r^2$ and on $r \, (\text{mod } 2N)$ (\cite{E-Z} Theorem 2.2).
 In (\ref{Fourier}), if the condition $4Nn-r^2\ge 0$ is deleted,
 we obtain a {\it nearly holomorphic Jacobi form}.
 \par Let $J_{k,N}^!$ be the space of
 nearly holomorphic Jacobi forms of weight $k$ and index $N$.
 Let $J_{*,*}^!$ be the ring of all nearly holomorphic Jacobi
 forms and $J_{ev,*}^!$ its even weight subring. Then $J_{ev,*}^!$
 is the free polynomial algebra over $M_*^!(\Gamma)={\mathbb
 C}[E_4,E_6,\Delta^{-1}]/(E_4^3-E_6^2=1728\Delta)$ on two
 generators $a=\tilde{\phi}_{-2,1}(\tau,z)\in
 {J}_{-2,1}^!$ and $b=\tilde{\phi}_{0,1}(\tau,z)\in
 {J}_{0,1}^!$ (for details, see \cite{E-Z}
 \S 9).
 Fix $k=2$ and $1\le N \le 6, \neq 4$.
 There are unique Jacobi forms $\phi_{D,N}\in J_{2,N}^!$
 having Fourier
 coefficients $c(n,r)=B(D, 4Nn-r^2)$ which depend only on the
 discriminant $r^2-4Nn$
 with $B(D,-D)=1$ and $B(D,d)=0$ if $d=4Nn-r^2<0,\neq -D$.
 The uniqueness of $\phi_{D,N}$ is obvious since
 the difference
 of any two functions satisfying the definition of $\phi_{D,N}$
 would be an element of $J_{2,N}$ (=the space of holomorphic
 Jacobi forms of weight $2$ and index $N$), which is of dimension zero
 by \cite{E-Z}, Theorem
 9.1 (2). For the existence,
 we need an additional condition on Fourier coefficients that
 $ B(D,0)=\begin{cases} -2, & \text{ \q if $D$ is a square } \\
                          0, & \text{ \q otherwise.} \end{cases}
 $
 \par\noindent
 The
 structure theorem then allows us to express $\phi_{D,N}$ as a linear
 combination of $a^i b^{N-i} (i=0,\dots,N)$ over $M_*^!(\Gamma)$.
 Define
 $$ g_{D,N}=q^{-D}+\sum_{d\ge 0} B(D,d) q^d. $$
 By the correspondence between Jacobi forms and half-integral
 forms (\cite{E-Z} Theorem 5.6),
 $g_{D,N}$ lies in the space $M_{3/2}^{+\cdots+}(N)^!$
 so that $f_{d,N}g_{D,N}$ defines a modular form of weight 2 for
 $\Gamma_0(4N)$. We write $f_{d,N}g_{D,N}=\sum_{n\in \mathbb Z}
 c_n q^n$.
 The ``plus" conditions imposed on $f_{d,N}$ and $g_{D,N}$
 force
 $(f_{d,N}g_{D,N})|_{U_{4N}}$ to be a modular form of weight 2
 on $SL_2(\mathbb Z)$.
 Here $U_{4N}$ is the operator sending $\sum_{n\in \mathbb Z}
 c_n q^n$ to $\sum_{n\in \mathbb Z}
 c_{4Nn} q^n$. In fact, if we consider
 $ h=\sum_{i\in (\mathbb Z/4N\mathbb Z)^\times} (f_{d,N}g_{D,N})
   \left( \frac{\tau + i}{4N} \right), $ then
 $h$ is invariant under the action of $\sm 0&-1\\1&0 \esm$ and
 has the Fourier development of the form
 $ \varphi (4N)\sum_{n\in \mathbb Z} e^{2\pi in/N} c_{4n}
 q^{n/N}$ since $c_n$ vanishes whenever $n\equiv 2 \mod 4$.
 $\sum_{i=0}^{N-1} h(\tau + i)=N\varphi (4N)
 (f_{d,N}g_{D,N})|_{U_{4N}}$ is then invariant under the action of
 $SL_2(\mathbb Z)$ with a pole only at $\infty$.
 Thus $(f_{d,N}g_{D,N})|_{U_{4N}}$ can be written as the
 derivative of some polynomial in $j$.
 By comparing the constant terms
 we get $A(D,d)=-B(D,d)$. This also shows the
  uniqueness of $f_{d,N}$.
 \par Through the article we adopt the following notations:
 \par\noindent $\bullet$ $T(m)$: generalized Hecke operator
 \par\noindent $\bullet$ $T_m$: Hecke operator acting on Jacobi
 forms or half-integral forms (\cite{E-Z} \S 4 and \S 5)
 \par\noindent $\bullet$ $\phi_D=\phi_{D,N}$
 \par\noindent $\bullet$ $g_D=g_{D,N}$
 \par\noindent $\bullet$ $f_d=f_{d,N}$
 \par\noindent $\bullet$ $\phi_D^{(p)}=\phi_{D,N^{(p)}}$
 \par\noindent $\bullet$ $g_D^{(p)}=g_{D,N^{(p)}}$
 \par\noindent $\bullet$ $B(d)=B(1,d)$

\section{Proof of Theorem \ref{Borcherds}}
 For each positive integer $m$ and prime $p$, we define
 $$J_m(d)=\sum_{Q\in {\mathcal Q}_{d,N}/\Gamma_0(N)^*}\frac{1}{w_Q} t_m(\alpha_Q)$$
 and
 $$J_m^{(p)}(d)
 =\sum_{Q\in {\mathcal Q}_{d,N^{(p)}}/\Gamma_0(N^{(p)})^*}\frac{1}{w_Q}
  t_m^{(p)}(\alpha_Q).$$
 First we need two lemmas.
 \begin{Lem} Let $p$ be a prime dividing $N$.
 For $i\ge 0$ and $m$ coprime to $p$,
 $$
 \phi_{p^{2i}m^2}^{(p)}|_{V_p}=p\phi_{p^{2i+2}m^2}+\phi_{p^{2i}m^2}.$$
 Here $V_p$ is the Hecke operator on Jacobi forms defined by the
 formula (2) in \cite{E-Z}.
 \label{AA}
 \end{Lem}
 \begin{proof}
 According to \cite{E-Z} Theorem 4.1, the operator $V_p$ maps
 $J_{2,N/p}^!$ to $J_{2,N}^!$.
 From the formula (7) in \cite{E-Z}, p.43, we find that
 $$
 \text{the coefficient of $q^n\zeta^r$ in
 }\phi_{p^{2i}m^2}^{(p)}|_{V_p}=
 \begin{cases} p, & \text{ if } 4Nn-r^2=-p^{2i+2}m^2 \\
                1, & \text{ if } 4Nn-r^2=-p^{2i}m^2 \\
                0, & \text{ if } 4Nn-r^2<0,
                    \neq -p^{2i+2}m^2,-p^{2i}m^2.\end{cases}$$
 From these observations and the uniqueness of $\phi_D$, the lemma
 immediately follows.
 \end{proof}
 \begin{Lem}
 Let $l$ be a positive integer coprime to $N$ and $d=4Nn-r^2$. Then
 \par\noindent
 (i) $J_l(d)= -\text{coefficient of $q^n\zeta^r$ in }\phi_1|_{T_l}$.
 \par\noindent
 (ii) $\phi_1|_{T_l}=\sum_{\nu | l} \nu \phi_{\nu^2}$.
 \label{BB}
 \end{Lem}
 \begin{proof} (i)
 Let $p$ be a prime divisor of $l$. Then
 \begin{align*}
 J_p(d) &=\sum_{Q\in {\mathcal Q}_{d,N}/\Gamma_0(N)^*}\frac{1}{w_Q} t_p(\alpha_Q)
   =\sum_{Q\in {\mathcal Q}_{d,N}/\Gamma_0(N)^*}\frac{1}{w_Q} t|_{T(p)}
   (\tau)|_{\tau=\alpha_Q} \\
  &=J_1(dp^2)+\left(\frac{-d}{p}\right)J_1(d)+pJ_1(d/p^2) \\
  & \text{ \q \q by a similar argument given in the proof of
  \cite{Zagier}
   Theorem 5-(ii) } \\
  &= -[B(dp^2)+\left(\frac{-d}{p}\right)B(d)+pB(d/p^2)] \text{ \q by
  \cite{Zagier} Theorem 8 } \\
  &= -B_p(d).
 \end{align*}
 Here $J_1(d/p^2)$ (resp. $B(d/p^2)$)
 is defined to be zero unless $d/p^2$ is an integer. And
 $B_p(d)$ denotes the coefficient of $q^d$ in $g_1|_{T_p}$,
 which is the same as the coefficient of $q^n\zeta^r$ in $\phi_1|_{T_p}$
 (\cite{E-Z} Theorems 4.5 and 5.4). Now let $p^s|| l$. Observe that
 $t|_{T(p^s)}=t_{p^{s-1}}|_{T(p)}-pt_{p^{s-2}}(\tau)$.
 Thus
 \begin{align*}
 J_{p^s}(d)
 &=J_{p^{s-1}}(dp^2)+\left(\frac{-d}{p}\right)J_{p^{s-1}}(d)+
   pJ_{p^{s-1}}(d/p^2)-pJ_{p^{s-2}}(d) \\
 &=-[B_{p^{s-1}}(dp^2)+\left(\frac{-d}{p}\right)B_{p^{s-1}}(d)+
   pB_{p^{s-1}}(d/p^2)]+pB_{p^{s-2}}(d) \text{ by induction on } s \\
 &=-\text{coefficient of $q^n\zeta^r$ in }
   [(\phi_1|_{T_{p^{s-1}}})|_{T_p}-p\phi_1|_{T_{p^{s-2}}}]\\
 &=-\text{coefficient of $q^n\zeta^r$ in }\phi_1|_{T_{p^{s}}}
  \text{ by Corollary 1 in \cite{E-Z} p.51.}
 \end{align*}
 Now write $l=l' p^s$ with $(l',p)=1$. Let $n(l)$ be the
 number of prime factors of $l$. We will use induction on
 $n(l)$.
 If $n(l)=1$, it returns to the previous case. Now
 \par\noindent
 $t|_{T(l)}=t|_{T(l')T(p^s)}=t_{l'}|_{T(p^s)}=
  t_{l'p^{s-1}}|_{T(p)}-pt_{l'p^{s-2}}$ which yields that
 \begin{align*}
 J_l(d)
 &=J_{l'p^{s-1}}(dp^2)+\left(\frac{-d}{p}\right)J_{l'p^{s-1}}(d)+
   pJ_{l'p^{s-1}}(d/p^2)-pJ_{l'p^{s-2}}(d) \\
 &=-\text{coefficient of $q^n\zeta^r$ in }\phi_1|_{T_{l}}
 \text{ by induction on } s.
 \end{align*}
 \par\noindent
 (ii) As before let $p$ be a prime dividing $l$ and $p^s||l$.
 First we will show that $\phi_1|_{T_{p^s}}=\sum_{i=0}^s p^i\phi_{p^{2i}}$.
 Let $s=1$. Then
 the coefficient of $q^d$ in $g_1|_{T_{p}}$ is
 $$B(dp^2)+\left(\frac{-d}{p}\right)B(d)+p B(d/p^2)=
  \begin{cases} 1, & \text{ if } d=-1 \\
                p, & \text{ if } d=-p^2 \\
                0, & \text{ if } d<0,\neq -1,-p^2. \end{cases}$$
 This implies $g_1|_{T_{p}}=pg_{p^2}+g_1$ and therefore
 $\phi_1|_{T_{p}}=p\phi_{p^2}+\phi_1$. Now let $s\ge 2$.
 Then
 \begin{align*}
 \phi_1|_{T_{p^s}}
 &=(\phi_1|_{T_{p^{s-1}}})|_{T_p} - p\phi_1|_{T_{p^{s-2}}} \\
 &=(\sum_{i=0}^{s-1} p^i\phi_{p^{2i}})|_{T_p}-
 p\sum_{i=0}^{s-2} p^i\phi_{p^{2i}}\text{ by induction on } s.
 \end{align*}
 For $i>0$, the coefficient of $q^d$ in $g_{p^{2i}}|_{T_{p}}$ is
 $$B(p^{2i},dp^2)+\left(\frac{-d}{p}\right)B(p^{2i},d)+p B(p^{2i},d/p^2)=
  \begin{cases} 1, & \text{ if } d=-p^{2i-2} \\
                p, & \text{ if } d=-p^{2i+2} \\
                0, & \text{ if } d<0,\neq -p^{2i-2},-p^{2i+2}.
                \end{cases}$$
 This shows that $\phi_{p^{2i}}|_{T_{p}}=
 \begin{cases} \phi_{p^{2i-2}}+p\phi_{p^{2i+2}}, & \text{ if } i>0
 \\
 \phi_1+p\phi_{p^2}, & \text{ if } i=0. \end{cases}$
  \q \q Thus
 \begin{align*}
 \phi_1|_{T_{p^s}}
 &=(\sum_{i=0}^{s-1} p^i\phi_{p^{2i}})|_{T_p}-
 p\sum_{i=0}^{s-2} p^i\phi_{p^{2i}}=\sum_{i=1}^{s-1}
 p^i(\phi_{p^{2i-2}}+p\phi_{p^{2i+2}})+\phi_1+p\phi_{p^2}
 -p\sum_{i=0}^{s-2} p^i\phi_{p^{2i}} \\
 &= \sum_{i=0}^s p^i\phi_{p^{2i}}.
 \end{align*}
 As in the proof of (i), write $l=l'p^s$ with $(l',p)=1$ and use
 induction on the number $n(l)$ of prime divisors of $l$.
 If $n(l)=1$, the assertion is clear.
 If $n(l)$ is greater than 1, then
 \begin{align*}
 \phi_1|_{T_{l}}
 &=\phi_1|_{T_{l'}T_{p^s}}=(\sum_{\nu|l'} \nu \phi_{\nu^2})|_{T_{p^s}}
 \text{ by induction on } n(l) \\
 &=\sum_{\nu|l} \nu \phi_{\nu^2} \text{ by induction on } s
  \text{ and applying the same argument as before.}
 \end{align*}
 \end{proof}
 We claim that for $d=4Nn-r^2$,
 \begin{equation}
  J_m(d)=-\text{coefficient of $q^n\zeta^r$ in }
  \sum_{u|m} 2^{s(u,N)}u\phi_{u^2}.
  \label{CC}
 \end{equation}
 Let $p$ be a prime dividing $N$.
 By \cite{Koike} Theorem 6.3 (2) (or \cite{Fer2} Proposition 2.6),
 the generalized Hecke operator $T(p)$ satisfies
 the following composition rule: for $k\ge 0$,
 $$
 T(p^k)\circ T(p)=T(p^{k+1})+pI_p\circ T(p^{k-1})
 $$ where
 $t_n|_{I_p}=t_n^{(p)}$ and
 $t_n$ is defined to be 0 if $n$ is not a rational integer.
 For $l$ coprime to $p$, we obtain
 \begin{align}
 t_{lp^{k+1}}
 &=t_l|_{T(p^{k+1})}=(t_l|_{T(p^{k})})|_{T(p)}-pt_l^{(p)}|_{T(p^{k-1})}
 \notag \\
 &=t_{lp^{k}}|_{T(p)}-pt_{lp^{k-1}}^{(p)}
  =t_{lp^{k}}^{(p)}(p\tau)+pt_{lp^{k}}|_{U_p}-pt_{lp^{k-1}}^{(p)}.
  \label{composition}
 \end{align}
 Meanwhile \cite{Koike} Theorem 3.1 Case I (or \cite{Fer2} Theorem
 3.7
 Case 1) provides the formula
 \begin{equation}
 pt_{lp^{k}}|_{U_p}+t_{lp^{k}}=t_{lp^{k}}^{(p)}+pt_{lp^{k-1}}^{(p)}.
 \label{compression}
 \end{equation}
 Combining (\ref{composition}) with (\ref{compression})
 we come up with $t_{lp^{k+1}}(\tau)=t_{lp^{k}}^{(p)}(p\tau)+t_{lp^{k}}^{(p)}(\tau)
  -t_{lp^{k}}(\tau)$ and therefore
 \begin{equation}
 \sum_{Q\in {\mathcal Q}_{d,N}/\Gamma_0(N)^*} t_{lp^{k+1}}(\alpha_Q)
 =\sum_{Q\in {\mathcal Q}_{d,N}/\Gamma_0(N)^*}
  (t_{lp^{k}}^{(p)}(p\tau)+t_{lp^{k}}^{(p)}(\tau))|_{\tau
  =\alpha_Q}
 -\sum_{Q\in {\mathcal Q}_{d,N}/\Gamma_0(N)^*}t_{lp^{k}}(\alpha_Q).
 \label{Hecke}
 \end{equation}
 The map which sends $[a,b,c]\in {\mathcal Q}_{d,N}$ to $[a/p,b,cp]\in
 {\mathcal Q}_{d,N/p}$ induces a bijection between ${\mathcal Q}_{d,N}/\Gamma_0(N)^*$
 and ${\mathcal Q}_{d,N/p}/\Gamma_0(N/p)^*$. And the natural map from
 ${\mathcal Q}_{d,N}/\Gamma_0(N)^*$ to ${\mathcal Q}_{d,N/p}/\Gamma_0(N/p)^*$ also gives
 a bijection. Thus (\ref{Hecke}) is rewritten as
 $$
 \sum_{Q\in {\mathcal Q}_{d,N}/\Gamma_0(N)^*} t_{lp^{k+1}}(\alpha_Q)
 =2\sum_{Q\in {\mathcal Q}_{d,N/p}/\Gamma_0(N/p)^*}
  t_{lp^{k}}^{(p)}(\alpha_Q)
 -\sum_{Q\in {\mathcal Q}_{d,N}/\Gamma_0(N)^*}t_{lp^{k}}(\alpha_Q),
 $$ which yields
 \begin{equation}
 J_{lp^{k+1}}(d)=2J_{lp^{k}}^{(p)}(d)-J_{lp^{k}}(d)
 \text{ \hspace{0.5cm} for }k\ge 0.
 \label{DD}
 \end{equation}
 We divide $N$ into two cases.
 \par\noindent Case I. $N=p=2$ or 3 or 5
 \par\noindent
 In (\ref{CC}) we write $m=l p^k$ with $(l,p)=1$.
 We use induction on $k$ to prove the claim.
 If $k=0$, the claim (\ref{CC}) follows from Lemma \ref{BB}.
 Now assume the claim for $k$.
 \begin{align*}
 &J_{lp^{k+1}}(d)
  =2J_{lp^{k}}^{(p)}(d)-J_{lp^{k}}(d) \\
 &=-\text{coefficient of $q^n\zeta^r$ in }
  [2(\phi_1^{(p)}|_{T_{lp^k}})|_{V_p}
  -(\sum_{i=1}^k\sum_{\nu |l}
  2\nu p^i\phi_{\nu^2p^{2i}}+\sum_{\nu|l}\nu\phi_{\nu^2})]\\
 & \text{ \hspace{1cm} by \cite{Zagier} Theorem 5-(ii) and induction
 hypothesis}
 \\
 &=-\text{coefficient of $q^n\zeta^r$ in }
  [2(\sum_{i=0}^k\sum_{\nu|l}\nu p^i\phi_{\nu^2p^{2i}}^{(p)})|_{V_p}
   -(\sum_{i=1}^k\sum_{\nu|l}
  2\nu p^i\phi_{\nu^2p^{2i}}+\sum_{\nu|l}\nu\phi_{\nu^2})]\\
 & \text{ \hspace{1cm} by \cite{Zagier} formula (19) and Theorem 5-(iii)
  }
 \\
 &=-\text{coefficient of $q^n\zeta^r$ in }
  [2\sum_{i=0}^k\sum_{\nu|l}\nu p^i(p\phi_{\nu^2p^{2i+2}}+\phi_{\nu^2p^{2i}})
  \\
 & \hspace{4.1cm} \q -(\sum_{i=1}^k\sum_{\nu|l}
  2\nu p^i\phi_{\nu^2p^{2i}}+\sum_{\nu|l}\nu\phi_{\nu^2})]
   \text{ \hspace{0.5cm} by Lemma \ref{AA} } \\
 &=-\text{coefficient of $q^n\zeta^r$ in }
  [\sum_{i=1}^{k+1}\sum_{\nu|l}
  2\nu p^i\phi_{\nu^2p^{2i}}+\sum_{\nu|l}\nu\phi_{\nu^2} ]
  \text{ \hspace{0.5cm} as desired.}
 \end{align*}

 \par\noindent Case II. $N=6$
 \par\noindent
 In (\ref{CC}) we write $m=l2^{k_1} 3^{k_2}$ with $(l,6)=1$ and ${k_1},{k_2}\ge0$.
 For simplicity, we
 put $\alpha(u)=2^{s(u,2)} u$ and $\beta(u)=2^{s(u,6)} u$.
 We will use induction on ${k_1}+{k_2}$.
 If ${k_1}+{k_2}=0$, the claim is immediate from Lemma \ref{BB}.
 Now assume ${k_1}+{k_2}\ge 1$, say ${k_2}\ge 1$.
 \begin{align*}
 & J_{l2^{k_1}3^{k_2}}(d)
  =2J_{l2^{{k_1}}3^{{k_2}-1}}^{(3)}(d)-J_{l2^{{k_1}}3^{{k_2}-1}}(d)
  \text{ \hspace{0.5cm} by (\ref{DD}) }\\
 &=-\text{coefficient of $q^n\zeta^r$ in }
  [2\sum_{i=0}^{k_1}\sum_{j=0}^{{k_2}-1}\sum_{\nu|l}
  \alpha(\nu 2^i3^j)\phi_{(\nu 2^i3^j)^2}^{(3)}|_{V_3}
  -\sum_{i=0}^{k_1}\sum_{j=0}^{{k_2}-1}\sum_{\nu|l}
  \beta(\nu 2^i3^j)\phi_{(\nu 2^i3^j)^2}] \\
 & \text{ \hspace{1cm} by the result in the case $N=2$
   and induction hypothesis } \\
 &=-\text{coefficient of $q^n\zeta^r$ in } \\
 & \q \q [2\sum_{i=0}^{k_1}\sum_{j=0}^{{k_2}-1}\sum_{\nu|l}
   \alpha(\nu 2^i3^j)(3\phi_{(\nu 2^i3^{j+1})^2}+\phi_{(\nu 2^i3^j)^2})
   -\sum_{i=0}^{k_1}\sum_{j=0}^{{k_2}-1}\sum_{\nu|l}
  \beta(\nu 2^i3^j)\phi_{(\nu 2^i3^j)^2}] \\
 & \text{ \hspace{1cm} by Lemma \ref{AA} } \\
 & =-\text{coefficient of $q^n\zeta^r$ in } \q
  \big[ \sum_{i=0}^{k_1}\sum_{\nu|l}
 2\alpha(\nu 2^i3^{{k_2}-1})\cdot 3\phi_{(\nu 2^i3^{k_2})^2}\\
 & \q \q +\sum_{i=0}^{k_1}\sum_{j=1}^{{k_2}-1}\sum_{\nu|l}
  [2\alpha(\nu 2^i3^{j-1})\cdot
  3+2\alpha(\nu 2^i3^{j})-\beta(\nu 2^i3^j)]\phi_{(\nu 2^i3^j)^2}\\
 & \q \q +\sum_{i=0}^{k_1}\sum_{\nu|l}[2\alpha(\nu 2^i)-\beta(\nu 2^i)]
 \phi_{(\nu 2^i)^2}\big]
 \\
 & =-\text{coefficient of $q^n\zeta^r$ in } \q
 \sum_{i=0}^{k_1}\sum_{j=0}^{{k_2}}\sum_{\nu|l}\beta(\nu 2^i3^j)\phi_{(\nu 2^i3^j)^2},
 \text{ \q as desired. }
 \end{align*}
 Let
 $z\in {\mathfrak H}$.
 Note that $\frac 1m t_m(z)$ can be viewed as the coefficient of $q^m$-term
 in $-\log q - \log (t(\tau)-t(z))$ (see \cite{Norton}). Thus
 $\log q^{-1}-\sum_{m > 0} \frac 1m t_m(z) q^m = \log (t(\tau)-t(z))$.
 Taking exponential on both sides, we get
 \begin{equation}
 q^{-1} \exp (-\sum_{m> 0} \frac 1m t_m(z) q^m)=t(\tau)-t(z).
 \label{product1}
 \end{equation}
 Define $B^*(u^2,d)=2^{s(u,N)} B(u^2,d)$. By the claim (\ref{CC}), we obtain
 \begin{equation}
 J_m(d)=-\sum_{u|m}uB^*(u^2,d).
 \label{EE}
 \end{equation}
  From (\ref{product1}) and (\ref{EE}) it follows that
 \begin{align*}
 {\mathcal H}_d(t(\tau))
 &=q^{-H(d)}\exp (-\sum_{m=1}^\infty J_m(d) q^m/m)
  =q^{-H(d)}\exp (\sum_{m=1}^\infty \sum_{u|m}uB^*(u^2,d) q^m/m) \\
 &=q^{-H(d)}\exp (\sum_{m=1}^\infty \sum_{u=1}^\infty
   uB^*(u^2,d) q^{mu}/(mu))\\
 &=q^{-H(d)}\exp (\sum_{u=1}^\infty (-B^*(u^2,d))\sum_{m=1}^\infty
   -(q^u)^m/m)\\
 &=q^{-H(d)}\exp (\sum_{u=1}^\infty \log (1-q^u)^{-B^*(u^2,d)})
  =q^{-H(d)}\prod_{u=1}^\infty (1-q^u)^{-B^*(u^2,d)}.
 \end{align*}
 Now the fact
 $A(D,d)=-B(D,d)$ completes the proof of our theorem.
 \begin{Rmk}
 If $N=4$, our proof does not apply since in this case the
 2-plicate $t^{(2)}$ of $t$ is the Hauptmodul for $\Gamma_0(2)$
 which is not $\Gamma_0(N)^*$-invariant for any $N$. In fact, we
 can numerically check that Theorem \ref{Borcherds} fails when
 $N=4$.
 \end{Rmk}
\section{Some recursion formulas}
\par
 Let $\delta$ be the denominator of $H(d)$. In the course of
 proving  Theorem \ref{Borcherds} we have seen that
 $${\mathcal H}_d(t(\tau))
 =q^{-H(d)}\prod_{m=1}^\infty\exp (-\sum_{u|m}uA^*(u^2,d) q^m/m).$$
 Observe that $(q^{H(d)}{\mathcal H}_d(t(\tau)))^\delta$ is of the
 form $1+\sum_{m=1}^\infty c(m)q^m$ with $c(m)\in \mathbb Z$.
 Then
 $$1+\sum_{m=1}^\infty c(m)q^m=\prod_{m=1}^\infty
   \exp (-\sum_{u|m}\delta uA^*(u^2,d) q^m/m).$$
 Put $V=\prod_{m=1}^\infty
   \exp (-\sum_{u|m}\delta uA^*(u^2,d) q^m/m).$ The differential
 identity $(\log V)'=V'/V$ (here $'$ denotes
   $q\frac{d}{dq}=\frac{1}{2\pi i}\frac{d}{d\tau})$
 leads to
 $$(-\sum_{m=1}^\infty \sum_{u|m}\delta uA^*(u^2,d) q^m)\cdot
   (1+\sum_{m=1}^\infty c(m)q^m )=\sum_{m=1}^\infty m c(m)q^m.$$
 Comparing the coefficients of $q^m$ on both sides we get
 $$ \sum_{u|m}\delta uA^*(u^2,d)+\sum_{1\le k <m} c(m-k)
    (\sum_{u|k}\delta uA^*(u^2,d))=-mc(m).$$
 Now we come up with the following
 recursion formula for $A^*(m^2,d)$: for $m\ge 1$,
 \begin{equation} A^*(m^2,d)=-\frac{1}{\delta} c(m) -\frac 1m
   [\sum_{1\le u<m \atop u|m} uA^*(u^2,d)+\sum_{1\le k <m}
   c(m-k)(\sum_{u|k} u A^*(u^2,d))].
 \label{recursion2}
 \end{equation}
 Thus all $A^*(m^2,d)$ can be computed from the values of $c(m)$.
 Likewise all $c(m)$ can be estimated recursively from the
 values of $A^*(m^2,d)$.
 \begin{Exm} \underline{$N=2, d=4$} \q
 Theorem \ref{Borcherds} yields the following product formula:
 \begin{equation}
  \left(t(\tau)-t\left(\frac{1+\sqrt{-1}}{2}\right)\right)^{1/2}
  =q^{-1/2}\prod_{u=1}^\infty (1-q^u)^{A^*(u^2,d)}.
  \label{FF}
 \end{equation}
 Here the Hauptmodul $t$ for $\Gamma_0(2)^*$ can be described by
 means of
 Dedekind $\eta$-functions, i.e.,
 \begin{align*}
 t(\tau)&=\left(\frac{\eta(\tau)}{\eta(2\tau)}\right)^{24}+24
    +\left(\frac{\eta(\tau)}{\eta(2\tau)}\right)^{24} \\
   &
   =q^{-1}+4372q+96256q^2+1240002q^3+10698752q^4+74428120q^5+\cdots,
 \end{align*}
 from which we obtain $t\left(\frac{1+\sqrt{-1}}{2}\right)=-104$.
 The identity (\ref{FF}) is then rewritten as
 $$
 1+104q+4372q^2+\cdots
 =\prod_{u=1}^\infty (1-q^u)^{2A^*(u^2,d)}
 =\prod_{m=1}^\infty
 \exp (-\sum_{u|m} 2uA^*(u^2,d) q^{m}/m).
 $$
 In (\ref{recursion2}) we take $\delta=2$, $c(1)=104$,
 $c(2)=4372$, $c(3)=96256$, etc. Then
 \begin{align*}
 A^*(1,4) & = -\frac 12 c(1) = -52, \\
 A^*(4,4) & = -\frac 12 c(2) -\frac 12[A^*(1,4)+c(1)A^*(1,4)]=544,
 \\
 A^*(9,4) & = -\frac 12 c(3) -\frac 13[A^*(1,4)+c(2)A^*(1,4)
             +c(1)A^*(1,4)+c(1)\cdot 2 \cdot A^*(4,4)]=-8244,\\
          & \q  \vdots
 \end{align*}
 \end{Exm}
 \begin{center}
 {\bf Appendix}
 \end{center}
 Let $f_0=\theta$.
 We found the initial $f_d$'s by expressing
 $[f_0, E_{12-2n}(4N\tau)]_{n}/\Delta(4N\tau)$ (if necessary,
 $[f_d, E_{12-2n}(4N\tau)]_n/\Delta(4N\tau))$ as linear combinations
 of them for $n=1,2,3,4$. Here
 $E_k$ is the normalized Eisenstein series of weight $k$, $\Delta$
 is the modular discriminant and
 $[ \q , \q ]_n$ denotes the ``Cohen
 bracket" (\cite{Cohen} \S 7 or \cite{Zagier94} \S 1).
 \par\noindent \underline{$N=2$}
 \par\noindent
 $f_4=q^{-4}-52q+272q^4+2600q^8-8244q^9
 +15300q^{12}+71552q^{16}-204800q^{17}+282880q^{20}+\cdots,$
 \par\noindent
 $f_7=q^{-7}-23q-2048q^4+45056q^8+252q^9
 -516096q^{12}+4145152q^{16}-1771q^{17}-26378240q^{20}+\cdots,$
 \par\noindent \underline{$N=3$}
 \par\noindent
 $f_3=q^{-3}-14q+40q^4-78q^9
 +168q^{12}-378q^{13}+688q^{16}+\cdots,$
 \par\noindent
 $f_8=q^{-8}-34q-188q^4+2430q^9
 +8262q^{12}-11968q^{13}-34936q^{16}+\cdots,$
 \par\noindent
 $f_{11}=q^{-11}+22q-552q^4-11178q^9
 +48600q^{12}+76175q^{13}-269744q^{16}
 +\cdots,$
 \par\noindent \underline{$N=5$}
 \par\noindent
 $f_4=q^{-4}-8q+q^4+10q^5+12q^9
 -62q^{16}+65q^{20}
 +\cdots,$
 \par\noindent
 $f_{11}=q^{-11}-12q-56q^4-45q^5+276q^9
 +672q^{16}+2520q^{20}
 +\cdots,$
 \par\noindent
 $f_{15}=q^{-15}-38q+112q^4-96q^5-988q^9
 +8512q^{16}+11856q^{20}
 +\cdots,$
 \par\noindent
 $f_{16}=q^{-16}-6q-132q^4+120q^5-1014q^9
 +3585q^{16}+17030q^{20}
 +\cdots,$
 \par\noindent
 $f_{19}=q^{-19}+20q+56q^4-210q^5-780q^9
 -23200q^{16}+46760q^{20}
 +\cdots,$
 \par\noindent \underline{$N=6$}
 \par\noindent
 $f_8=q^{-8}-10q-12q^4+54q^9
 +54q^{12}-88q^{16}+\cdots,$
 \par\noindent
 $f_{12}=q^{-12}-28q+26q^4-156q^9
 +168q^{12}+728q^{16}
 +\cdots,$
 \par\noindent
 $f_{15}=q^{-15}-10q-64q^4+3q^9
 -320q^{12}+1664q^{16}
 +\cdots,$
 \par\noindent
 $f_{20}=q^{-20}+12q-64q^4-756q^9
 +945q^{12}-2912q^{16}
 +\cdots,$
 \par\noindent
 $f_{23}=q^{-23}-13q+64q^4-27q^9
 -1728q^{12}-5760q^{16}
 +\cdots,$
 \par\noindent
 For the remaining $f_d(\tau)$ we inductively obtain them by multiplying
 $f_{d-4N}(\tau)$ by $j(4N\tau)$ to get a ``plus" form of weight
 $1/2$ with leading coefficient $q^{-d}$ and then subtracting a
 suitable linear combination of $f_{d'}(\tau)$ with $0\le d' <d$.
 \vspace{0.2cm}

 \begin{center}
{\bf Acknowledgment}
\end{center}
 I am grateful to Professor Don Zagier for introducing me to this
 subject. I would also like to take an opportunity to
 thank Professor Richard E. Borcherds, Professor Jan H.
 Bruinier and Professor Ja Kyung Koo for their kind and valuable comments.
 
\end{document}